\documentclass{article}
\usepackage[OT1,T1]{fontenc}
\usepackage{fancyhdr}
\usepackage[counterclockwise]{rotating}
\usepackage{amsmath}
\usepackage{latexsym,amssymb}
\usepackage{amsfonts}
\usepackage{hyperref}
\usepackage{url}
\usepackage{amsthm}
\usepackage{mathrsfs}
\usepackage{ifthen}

\usepackage{xy}
\usepackage{amscd}

\newtheorem{theo}{Theorem}[section]

\newtheorem{corol}{Corollary}[section]
\newtheorem{lem}{Lemma}[section]

\theoremstyle{change}

\newtheorem{rem}{Remark}[section]

\xyoption{all} \xyoption{matrix}

\newcommand{\Hu}{H^{1}\left(\Omega\right)}

\newcommand{\Huz}{H^{1}_{0}\left(\Omega\right)}

\newcommand{\Ld}{L^{2}\left(\Omega\right)}
\newcommand{\U} {H^2(\Omega)\cap H^1_0(\Omega)}%
\newcommand{\Hd}{H^{2}\left(\Omega\right)}
\newcommand{\m}{\text{meas}}
\newcommand{\ess}{\text{ess}}

\DeclareGraphicsExtensions{eps}

\def\keywords#1{\parbox{\hsize}{
\tolerance=200 \parindent=10pt {\bf Key words:} \small
#1.} \vskip1pt}%

\def\subclass#1{\parbox{\hsize}{
\tolerance=200 \parindent=10pt {\small\bf AMS Classification
(2000):} \small #1} \vskip1pt\par\noindent}

\begin{document}
%--------------------------------------------------------------------------------------------------------
\title{\textbf {Optimal control of unilateral obstacle problem with a source term}}

%---------------------------------------------------------------------------------------------------------

\author{\textsc{Radouen Ghanem}  \\
{\footnotesize  {UMR 6628-MAPMO, Fédération Denis Poisson,
Université d'Orléans}}\\ {\footnotesize BP. 6759, F-45067 Orléans
Cedex 2}\\ {\footnotesize{ \texttt{radouen.ghanem@univ-orleans.fr}
}} }
\date{}
\maketitle

\begin{abstract}

We consider an optimal control problem for the obstacle problem
with an elliptic variational inequality. The obstacle function
which is the control function is assumed in $H^{2}$. We use an
approximate technique to introduce a family of problems governed
by variational equations. We prove optimal solutions existence and
give necessary optimality conditions.
\end{abstract}\bigskip

%\subclass{??}%

\subclass{\small 35R35, 49J40, 49J20.}%

%-----------------------------------------------------------------------
\section{Introduction}
%------------------------------------------------------------------------

\indent The study of variational inequalities and free boundary
problems finds application in a variety of disciplines including
physics, engineering, and economies as well as potential theory
and geometry. In the past years, the optimal control of
variational inequalities has been studied by many authors with
different formulations. For example optimal control problems for
obstacle problems (where the obstacle is a given (fixed) function)
were considered with the control variables in the variational
inequality. Roughly speaking, the control is different from the
obstacle, see for example works by \cite{Barbu3}, \cite{Bergou3},
\cite{Mign}, \cite{MignPuel}  and the references therein.

Here we deal with the obstacle as the control function. This kind
of problem appears in shape optimization for example. It may
concern a dam optimal shape. The obstacle gives the form to be
designed such that the pressure of the fluid inside the dam is
close to a desired value. This is equivalent in some sense to
controlling the free boundary \cite{Friedman}.

The main difficulty of this type of problem comes from the fact
that the mapping between the control and the state
(control-to-state operator) is not differentiable but only
Lipschitz-continuous and so it is not easy to get first order
optimality conditions.

These problems have been considered from the theoretical and/or
numerical points of view by many authors (see for example Adams
and Lenhart \cite{Adam2}, Ito and Kunisch \cite{Kunish}). They
have used either an approximation of the variational inequality by
penalization-regularization or a complementarity constraint
formulation. Adams et Lenhart \cite{Adam2} consider optimal
control problem governed by a linear elliptic variational
inequality without source terms. The main result is that any
optimal pair must satisfy "state = obstacle". Adams and Lenhart
\cite{Adam} treat control of $H^{1}-$ obstacle, and convergence
results in \cite{Adam}, \cite{Adam2} are given under implicit
monotonicity assumptions.

Ito and Kunisch \cite{Kunish} consider the optimal control problem
to minimize a functional involving the $H^{1}$ norm of the
obstacle, subject to a variational inequality of the type
$y\in{\rm argmin}\{a(z)-\langle f,z\rangle| z\le\psi\}$ in a
Hilbert lattice $H$. Under appropriate conditions, they show that
the variational inequality can be expressed by the system
$Ay+\lambda=f,
\lambda\mathrel{\mathop:}=\max(0,\lambda+c(y-\psi))$. Smoothing
the max-operation, this system is approximated by a semilinear
elliptic equation containing only smooth expressions. Passing to
the limit, the optimality system of the associated differentiable
optimal control problem is used to derive an optimality system of
the original nonsmooth control problem with only
$H^{1}$-regularity for the obstacle.

Bergounioux and Lenhart \cite{Bergou1}, \cite{Bergou2} have
studied obstacle optimal control for semilinear and bilateral
obstacle problem, where the admissible controls (obstacles) are
$H^{2}$-bounded and the convergence results are given with a
compactness assumption. Yuquan and Chen \cite{Yuquan}, consider an
obstacle control problem in a elliptic variational inequality
without source terms. We can see also quote the paper of Lou
\cite{Lou} for more generalized regularity results. In this paper
we consider an optimal control problem: we seek an optimal pair of
optimal solution (state, control), when the state is close to a
desired target profile and satisfies an unilateral variational
inequality with a source term, and the control function is the
lower obstacle. Convergence results are proved with compactness
techniques.

The new feature in this paper is the regularity on the control
function (obstacle) and an optimality conditions system more
complete than the one given in \cite{Yuquan}.

Let us give the outline of the paper. Next section, is devoted to
the formulation of the optimal control problem, we give
assumptions for the state equation, and give preliminaries
results. In section 3 we study the variational inequality, give
control-to-state operator properties and assert an existence
result for optimal solution. The last section is devoted to
optimality condition system.

\setcounter{equation}{0}
%----------------------------------------------------
\section{Optimal control problem}
%-------------------------------------------------------

Let $\Omega $ be an open bounded set in $\mathbb{R}%
^{n}$ $( n\leq 3) $, with lipschitz boundary $\partial \Omega$. We
adopt the standard notation $H^{m}( \Omega) $ for the Sobolev
space of order $m$ in $\Omega $ with norm $\left\Vert \cdot
\right\Vert _{H^{m}\left(\Omega\right)}$, where

\begin{equation*}
H^{m}( \Omega ) \mathrel{\mathop:} =\left\{ v\mid v\in
L^{2}(\Omega) ,\partial ^{q}v\in L^{2}( \Omega) ~\forall
q,\left\vert q\right\vert \leqslant m\right\},
\end{equation*}

 and

\begin{equation*}
H_{0}^{m}(\Omega) \mathrel{\mathop:} =\left\{ v\mid v\in
H^{m}(\Omega) ,\left. \frac{\partial ^{k}v}{\partial \eta
^{k}}\right\vert _{\partial \Omega }=0,~0 \leq k\leq m-1\right\},
\end{equation*}

defined as the closure of $\mathscr{D}\left(\Omega\right)$ in the
space $H^{m}\left(\Omega\right)$, where
$\mathscr{D}\left(\Omega\right)$, the space of
$\mathcal{C}^{\infty}\left(\Omega\right)$-functions, with compact
support in $\Omega$ (see for example \cite{Adams}). We shall
denote by $\left\Vert \cdot \right\Vert _{V}$ , the Banach space
$V$ norm, and  $\left\Vert \cdot \right\Vert
_{L^{p}\left(\Omega\right)}$ the $p-\text{summable}$ functions
$u:\Omega \rightarrow \mathbb{R}$ endowed with the norm
$\displaystyle
\|u\|_{L^{p}\left(\Omega\right)}\mathrel{\mathop:}=\left(\int_{\Omega}|u\left(x\right)|^{p}dx\right)^{1/p}$
for $1\leq p < \infty$ and
$\|u\|_{L^{\infty}\left(\Omega\right)}\mathrel{\mathop:}=\underset{x\in
\Omega}{\ess\sup}|u\left(x\right)|$ for $p=\infty$. In the same
way, $\left\langle \cdot ,\cdot \right\rangle $ denotes the
duality product between $H^{-1}( \Omega) $ and $H_{0}^{1}(\Omega)
$, and $( \cdot ,\cdot)$ the $L^{p}(\Omega) $ inner product. It is
well known that $H_0^1(\Omega)\hookrightarrow L^2(\Omega)
\hookrightarrow
H^{-1}(\Omega)$ with compact and dense injection. We consider the bilinear form $a( \cdot ,\cdot) $ defined on $%
H_{0}^{1}(\Omega) \times H_{0}^{1}(\Omega) $ by

\begin{equation}
a(u,v) \mathrel{\mathop:} = \sum_{i,j=1}^n \int_\Omega a_{ij} {\frac{\partial u }{\partial x_i }%
} {\frac{\partial v }{\partial x_j }}~d{x} + \sum_{i=1}^n \int_\Omega a_i {%
\frac{\partial u }{\partial x_i }} v ~d{x} +\int_\Omega a_{0} u\,
v ~d{x}, \label{bilfora}
\end{equation}

%where $a_{ij},a_{i},a_{0}$ belong to $L^{\infty }( \Omega)$.

where

\begin{equation*}
 \left\{
\begin{tabular}{l}
$a_{0},a_{i},a_{ij}\in L^{\infty }\left( \Omega \right),$ \\
\\
$\sum\limits_{i,j=1}^{n}a_{ij} \theta_{i}\theta _{j}\geq
m\sum\limits_{i=0}^{n}\theta _{i}^{2},\quad m>0,\text{\quad a.e.
in }\Omega ,\quad
\forall \theta \in \mathbb{R}^{n}.$ \\
\end{tabular}%
\right. \eqno{({\mathbf{H}})}
\end{equation*}

Moreover, we suppose that $a_{ij}\in \mathcal{\
C}^{0,1}(\bar{\Omega})$\ (the space of lipschitz continuous
functions in $\Omega $, where $\bar{\Omega}$ is the closure of
$\Omega$) and that $a_{0}$ is nonnegative to ensure a good
regularity of the solution (see for example \cite{LionMage}). We
suppose that the bilinear form $a( \cdot ,\cdot)$ is continuous on
$H_{0}^{1}( \Omega) \times H_{0}^{1}(\Omega)$

\begin{equation}
\exists\,M>0,\forall \varphi ,\psi \in H_{0}^{1}(\Omega)
,\left\vert a( \varphi ,\psi) \right\vert \leq M\left\Vert \varphi
\right\Vert _{H_{0}^{1}(\Omega) }\left\Vert \psi \right\Vert
_{H_{0}^{1}( \Omega) }, \label{conbilfora}
\end{equation}

and coercive on $H_{0}^{1}( \Omega) \times H_{0}^{1}(\Omega)$

\begin{equation}
\exists\,m>0,\forall \varphi \in H_{0}^{1}(\Omega) ,a(\varphi
,\varphi) \geq m\left\Vert \varphi \right\Vert
_{H_{0}^{1}(\Omega)}^2. \label{coebilfora}
\end{equation}

We call $A\in \mathcal{\ L}({H_{o}^{1}(\Omega )},H^{-1}(\Omega ))$
the linear (elliptic) operator associated to $a(\cdot ,\cdot) $
such that $\left\langle Au,v\right\rangle \mathrel{\mathop:}=
a(u,v)$. We note that the coercivity assumption \eqref{coebilfora}
on $a$ implies that

\begin{equation*}
\forall ~\varphi \in H_{0}^{1}(\Omega) ,~\left\langle A\varphi
,\varphi \right\rangle \geq m \left\Vert \varphi \right\Vert
_{H_{0}^{1}(\Omega) }^2.
\end{equation*}

For any $\varphi \in \Huz$, we define

\begin{equation*}
\mathcal{K}(\varphi)\mathrel{\mathop:} =\left\{ y \in
H_{0}^{1}(\Omega ) \mid y \geq \varphi~\mbox{a.e. in}~\Omega
\right\}, \label{K}
\end{equation*}

and consider the following variational inequality

\begin{equation}
a\left( y,\,v-y\right) \geq \left(f,\, v-y\right),\quad \forall v
\in \mathcal{K} (\varphi), \label{var_eq}
\end{equation}

where $f$ belongs to $\Ld$ as a source term. In addition $c$ or
$C$ denotes a general positive constant independent of $\delta$.

\begin{theo}
Under the hypothesis  \eqref{conbilfora} and \eqref{coebilfora},
for any $f$ $\in$ $L^{2}\left(\Omega\right)$ and $\varphi$ $\in$
$H_{0}^{1}\left( \Omega \right)$, the variational inequality
\eqref{var_eq}, has a unique solution
 $y$ in $\mathcal{K}\left(\varphi\right)$. In addition if $\varphi$ belongs to $H^{2}\left( \Omega \right)$, the solution
$y$ belongs to  $H^2(\Omega)\cap H^1_0(\Omega)$.
\end{theo}

\begin{proof}
See \cite{Brezi}.
\end{proof}

From now we define the operator $\mathcal{T}$ (control-to-state)
from $H^2(\Omega)\cap H^1_0(\Omega)$ to $H^2(\Omega)\cap
H^1_0(\Omega)$, such that
$y\mathrel{\mathop:}=\mathcal{T}\left(\varphi\right)$ is the
unique solution to the variational inequality \eqref{var_eq}.

Now, we consider the optimal control problem $\left(
\mathcal{P}\right)$, defined as follows

\begin{equation}
\min \left\{ J(\varphi
)\mathrel{\mathop:}=\tfrac{1}{2}\int\nolimits_{\Omega }\left(
\mathcal{T}\left( \varphi \right) -z\right) ^{2}d{
x}+\tfrac{\nu}{2}\left( \int\nolimits_{\Omega }\left( \Delta
\varphi \right) ^{2}d{x}\right),~\varphi \in
\mathcal{U}_{ad}\right\}, \tag{$\mathcal{P}$}
\end{equation}

where $\nu$ is a given positive constant, $z$ $\in$ $\Ld$ and
$\mathcal{U}_{ad}$ (the set of admissible control) is a closed
convex subset of $\U$: we seek an obstacle (optimal control)
$\varphi^{\ast}$ in $\mathcal{U}_{ad}$, such that the
corresponding state is close to a target profile $z$. In the
sequel we set $\mathcal{U} \mathrel{\mathop:}=\U$.

\setcounter{equation}{0}

%-------------------------------------------------------------------------------------------
\section{Approximation of problem ($\mathcal{P}$)}
\subsection{Approximation of operator $\mathcal {T}$}
%--------------------------------------------------------------------------------------------

The obstacle problem \eqref{var_eq} can be equivalently written as
follows

\begin{eqnarray}
\label{ine_sub_diff}%
Ay+\partial I_{\mathcal{K\left(\varphi\right)}}\left(y\right) \ni
f \quad \mbox{in}~\Omega, \quad y =0~\mbox{on}~\partial \Omega,
\end{eqnarray}

where

\begin{equation*}
\partial I_{\mathcal{K\left({\varphi}\right)}}\left(y\right)= \partial I_{\mathcal{K}^{+}\left(y - \varphi\right)}\left(y\right)
\mathrel{\mathop:}=\left\{ v\in L^{2}\left( \Omega \right) \mid v
\in \beta_{o} \left( y-\varphi  \right),\,\text{a.e. in}~\Omega
\right\},
\end{equation*}

and

\begin{equation*}
\mathcal{K}^{+}\mathrel{\mathop:} =\left\{ y \in H_{0}^{1}(\Omega
) \mid y \geq 0,~\mbox{a.e. in}~\Omega \right\}, \label{K}
\end{equation*}

and $\beta_{o} :\mathbb{R\longrightarrow }2^{\mathbb{R}}$ is the
maximal monotone (multivalued) graph,

\begin{equation*}
\beta_{o}(r)\mathrel{\mathop:}=
\begin{cases}
0 & \mathrm{if\ }r \geq 0 \\
\mathbb{R}^{-} & \mathrm{if\ }r=0 \\
\emptyset & \mathrm{if\ } r<0.%
\end{cases}%
\end{equation*}

Equation \eqref{ine_sub_diff}, can be approximated by the
following smooth  semilinear equation

\begin{equation}
Ay+\beta_{\delta} \left( y-\varphi \right) =f\quad
\mbox{in}~\Omega ,\quad y=0~\mbox{on}~\partial \Omega,
\label{semi_line_equ}
\end{equation}

where $\beta_{\delta}$  is an approximation of $\beta_{o}$. One
possible approximation of $\beta_{o}$ si given as follow

\begin{equation*}
\label{beta} \beta_{\delta}(r)\mathrel{\mathop:}=
\tfrac{1}{\delta}
\begin{cases}
0 & \mathrm{if\ }r \geq 0 \\
-r^{2} & \mathrm{if\ }r\in \left[ -\frac{1}{2},0\right] \\
r+\frac{1}{4} & \mathrm{if\ } r \leq -\frac{1}{2},%
\end{cases}%
\end{equation*}

Where $\delta > 0$ and we note that
$r_{\delta}^{-}\mathrel{\mathop:}=\tfrac{1}{\delta}\min\{0,r\}
\leq \beta_{\delta}\left(r\right) \leq 0$ and $\beta \in
\mathscr{C}^{\infty}\left( \mathbb{R}\right)$ and
$\beta_{\delta}^{\prime}$ is given by

\begin{equation*}
\beta_{\delta} ^{\prime }(r)\mathrel{\mathop:}= \tfrac{1}{\delta}
\begin{cases}
0 & \mathrm{if\ }r\geq 0 \\
-2r & \mathrm{if\ }r\in \left[ -\frac{1}{2},0\right] \\
1 & \mathrm{if\ }r\leq -\frac{1}{2}.%
\end{cases}%
\label{betaprime}
\end{equation*}

As $\beta_{\delta}(\cdot -\varphi)$ is nondecreasing, it is well
known (see for example \cite{Gilbarg}), that the boundary value
problem (\ref{semi_line_equ}) admits a unique solution $y^{\delta
}$ in $\U$ for a fixed $\varphi$ in $\U$ and $f$ in $\Ld$. In the
sequel, we set $y^{\delta }\mathrel{\mathop:}=\mathcal{T}^{\delta
}\left( \varphi \right)$. We recall the following continuity
results \cite{Bergou2}

\begin{theo}
\label{yhyleqphiphih} For any pair $({y}_{i},{\varphi}_{i})$ in
$\mathcal{U} \times \mathcal{U}$, that satisfies
\eqref{semi_line_equ} where $i=1,\,2$. We get
\begin{equation*}
\label{est_y_phi} \left\|{y}_{2}-{y}_{1}\right\|_{\Hu}\leq
L_{\delta}\left\| {\varphi}_{2}-{\varphi}_{1}\right\|_{\Ld},
\end{equation*}
where $L_{\delta}\mathrel{\mathop:}=\max \left\{
1,\tfrac{2}{m\delta}\right \}$ and $m$ is the coercivity constant
of $a$.
\end{theo}

\begin{proof}

From \eqref{statequavar}, we obtain

\begin{equation*}
a\left( y_{2}-y_{1},v\right) +\left( \beta _{\delta }\left(
y_{2}-\varphi_{2}\right) -\beta _{\delta }\left(
y_{1}-\varphi _{2}\right) ,v\right)=0,\quad\forall v \in \mathcal{U}.%
\end{equation*}

with $v=y_{2}-y_{1}$, we write

\begin{equation*}
a\left( y_{2}-y_{1},y_{2}-y_{1}\right) +\left( \beta _{\delta
}\left( y_{2}-\varphi_{2}\right) -\beta _{\delta }\left(
y_{1}-\varphi _{2}\right) ,y_{2}-y_{1}\right)=0.%
\label{equstatdif}
\end{equation*}

Since  $\beta_{\delta}$ is nondecreasing, by the hypothesis
\eqref{conbilfora} and \eqref{coebilfora}, we deduce

\begin{equation*}
\left\| y_{2}-y_{1}\right\|_{\Hu}\leq L_{\delta} \left\|
\varphi_{2}-\varphi_{1}\right\|_{\Ld},
\end{equation*}

where $L_{\delta}\mathrel{\mathop:}=\max \left\{
1,\tfrac{2}{m\,\delta}\right \}$.
\end{proof}

\begin{theo}
\label{ydcony} Let $\varphi ^{\delta }$ in $\U$ be a strongly
convergent sequence in $\Huz$ to some $\varphi$ as $\delta $ tends
to $0$. Then the sequence
$y^{\delta}\mathrel{\mathop:}=\mathcal{T}^{\delta }\left(
\varphi^{\delta} \right)$ strongly converges to
$y\mathrel{\mathop:}=\mathcal{T}\left( \varphi \right)$ in $\Huz$.
\end{theo}

\begin{proof}
For every $\varphi ^{\delta }$ in $\U$, we set $y^{\delta
}\mathrel{\mathop:}=\mathcal{T}^{\delta }\left( \varphi ^{\delta
}\right)$, then for any $y^{\delta}$ in $\Huz$ the equation
(\ref{semi_line_equ}) is equivalent to

\begin{equation}
a\left( y^{\delta },v\right) +\left( \beta_{\delta } \left(
y^{\delta }-\varphi ^{\delta }\right) ,v\right)=\left( f,v\right),
\quad \forall v\in H_{0}^{1}\left( \Omega \right) .
\label{semi_line_equ_var}
\end{equation}

In the equation (\ref{semi_line_equ_var}), we choose $v=\varphi
^{\delta }-y^{\delta }$, then we get

\begin{equation*}
a\left( y^{\delta },\varphi ^{\delta }-y^{\delta }\right)
+\int\nolimits_{\Omega }\beta_{\delta } \left( y^{\delta }-\varphi
^{\delta }\right) \left( \varphi ^{\delta }-y^{\delta }\right)
d{x}=\int\nolimits_{\Omega }f\left( \varphi ^{\delta }-y^{\delta
}\right) d{x}.
\end{equation*}

We know by the definition of $\beta_{\delta}$, thats if $y^{\delta
}\left( x\right) -\varphi ^{\delta }\left(x\right) \geq 0$,  we
have

\begin{equation*}
\beta_{\delta} \left( y^{\delta }\left(x\right) -\varphi ^{\delta
}\left( x\right) \right)=0,
\end{equation*}

otherwise, we get $\beta_{\delta} \left( y^{\delta }\left(x\right)
-\varphi ^{\delta }\left( x\right) \right) \leq 0$. Then we deduce
that in all cases, we have

\begin{equation*}
\beta_{\delta} \left( y^{\delta }-\varphi ^{\delta }\right) \left(
y^{\delta }-\varphi ^{\delta }\right) \geq 0\quad \mbox{a.e.
in}~\Omega ,
\end{equation*}

that yields

\begin{equation*}
a\left( y^{\delta },y^{\delta }\right) \leq a\left( y^{\delta
},\varphi ^{\delta }\right) +\left( f,y^{\delta }-\varphi ^{\delta
}\right),
\end{equation*}

with the hypothesis \eqref{conbilfora} and \eqref{coebilfora}, we
deduce (estimation regularity)

\begin{equation}
\left\Vert y^{\delta }\right\Vert_{\Hu}\leq C \left\Vert \varphi
^{\delta }\right\Vert_{\Hu}, \label{born_y_delta}
\end{equation}

where $C$ is a constant only depending on $f$ and $a$. We know
that if $\varphi ^{\delta }$ is strongly convergent in
$H_{0}^{1}\left( \Omega \right)$ then $\varphi ^{\delta }$ is
bounded in $H_{0}^{1}\left( \Omega \right)$, and by
(\ref{born_y_delta}) we deduce that $y^{\delta }$ is convergent to
some $y$ as $\delta$ tends to $0$ weakly in $H_{0}^{1}\left(
\Omega \right) $ and strongly in $L^{2}\left( \Omega \right)$.

Let $v$ in $\mathcal{K}\left( \varphi \right)$, and choose
$v^{\delta }=\max \left( v,\varphi ^{\delta }\right) $. We have
$v^{\delta }$ in  $\mathcal{K}\left( \varphi ^{\delta }\right) $
and that $v^{\delta }$ is convergent to $v$ strongly in
$H_{0}^{1}\left( \Omega \right)$. Equation
(\ref{semi_line_equ_var}) with $v=v^{\delta }-y^{\delta }$ gives

\begin{equation*}
a\left( y^{\delta },v^{\delta }-y^{\delta }\right)
+\int\nolimits_{\Omega }\beta_{\delta } \left( y^{\delta }-\varphi
^{\delta }\right) \left( v^{\delta }-y^{\delta }\right) d{
x}=\int\nolimits_{\Omega }f\left( v^{\delta }-y^{\delta }\right)
d{x}.
\end{equation*}

\begin{itemize}
\item If $y^{\delta }\leq \varphi ^{\delta }$, therefore
$\beta_{\delta} \left( y^{\delta }-\varphi ^{\delta }\right) <0$
and $\left( v^{\delta }-y^{\delta }\right) \geq 0$, we deduces
that $ \beta_{\delta} \left( y^{\delta }-\varphi ^{\delta }\right)
\left( v^{\delta }-y^{\delta }\right) \leq 0$.

\item If $y^{\delta }\geq \varphi^{\delta} $, then $y^{\delta
}-\varphi ^{\delta }\geq 0$, therefore $\beta_{\delta} \left(
y^{\delta }-\varphi ^{\delta }\right)=0$.
\end{itemize}

So we deduce that in all cases we have $\beta_{\delta} \left(
y^{\delta }-\varphi ^{\delta }\right) \left( v^{\delta }-y^{\delta
}\right) \leq 0$, and we get

\begin{equation*}
a\left( y^{\delta},y^{\delta }\right) \leq  a\left( y^{\delta
},v^{\delta }\right) -\left( f,v^{\delta }-y^{\delta }\right).
\end{equation*}

Passing to the limit and using the lower semi-continuity of $a$
gives

\begin{equation*}
a\left( y,y\right) \leq \underset{\delta \rightarrow 0}{\lim \inf~
}a\left( y^{\delta },y^{\delta }\right) \leq \underset{\delta
\rightarrow 0}{\lim \inf }~a\left( y^{\delta },v^{\delta }\right)
-\left( f,v^{\delta }-y^{\delta }\right)=a\left( y,v\right)
-\left( f,v-y\right),
\end{equation*}

and

\begin{equation*}
a\left( y,v-y\right) \geq \left( f,v-y\right),\quad \forall v \in
\mathcal{K} \left( \varphi \right).
\end{equation*}

It remains to prove that $y^{\delta}$ tends to $y$, strongly in
$H_{0}^{1}\left( \Omega \right) $. By using the fact that
$w^{\delta }=\max \left( y,\varphi ^{\delta }\right) $ converge to
$y$ strongly in $H_{0}^{1}\left( \Omega \right) $ it is sufficient
to prove that $w^{\delta }-y^{\delta }$ converge to $0$ strongly
in $H_{0}^{1}\left( \Omega \right)$. From equation
(\ref{semi_line_equ_var}) we get

\begin{equation*}
\begin{split}
a\left( w^{\delta }-y^{\delta },w^{\delta }-y^{\delta }\right) & =
 a\left( w^{\delta },w^{\delta }-y^{\delta }\right) -a\left(
y^{\delta },w^{\delta }-y^{\delta }\right) \\
& =  a\left( w^{\delta },w^{\delta }-y^{\delta }\right)
+\int\nolimits_{\Omega }\beta_{\delta } \left( y^{\delta }-\varphi
^{\delta }\right) \left( w^{\delta }-y^{\delta }\right) d{x} -
\int\nolimits_{\Omega }f\left( w^{\delta }-y^{\delta }\right)
d{x}.
 \end{split}
\end{equation*}

As previously we deduce that

\begin{equation*}
a\left( w^{\delta }-y^{\delta },w^{\delta }-y^{\delta }\right)
\leq  a\left( w^{\delta },w^{\delta }-y^{\delta }\right) -\left(
f,w^{\delta }-y^{\delta }\right),
\end{equation*}

\noindent from the hypothesis \eqref{coebilfora}, we get
\begin{equation*}
m \left\Vert w^{\delta }-y^{\delta }\right\Vert_{\Hu}^{2}  \leq
 a\left( w^{\delta },w^{\delta }-y^{\delta }\right) -\left(
f,w^{\delta }-y^{\delta }\right).
\end{equation*}

\end{proof}

As a consequence of the previous theorem, we obtain the following
corollaries

\begin{corol}

For any $\varphi^{\delta}$ in $\Hd \cap \Huz$,
$y^{\delta}\mathrel{\mathop:}=\mathcal{T}^{\delta}\left(
\varphi^{\delta} \right)$ belongs to $\Hd \cap \Huz$.
\end{corol}

\begin{proof}
Since $\beta_{\delta}\left(y^{\delta}-\varphi^{\delta}\right)$
 and $f$ belongs to $\Ld$, then $Ay^{\delta} \in \Ld$ and $y^{\delta}\in
 \Hd$.
\end{proof}

\begin{corol}
\label{ydconfy} For any $\varphi$ in $\mathcal{U}_{ad}$, the
sequence $y^{\delta}\mathrel{\mathop:}=\mathcal{T}^{\delta }\left(
\varphi \right) $, converges to
$y\mathrel{\mathop:}=\mathcal{T}\left( \varphi \right) $, strongly
in $H_{0}^{1}\left( \Omega \right)$.
\end{corol}

\begin{corol}
\label{cortocon} There exists a constant $C$ depending only on $f$
and $a$, such that for any $\varphi $ in $\mathcal{U}$, we get
\begin{equation*}
\left\| \mathcal{T}\left( \varphi \right) \right\|_{\Huz}\leq
C(a,f)\left\| \varphi \right\|_{\Huz}.
\end{equation*}
\end{corol}

\begin{proof}
We choose $\varphi ^{\delta}=\varphi $, and $y^{\delta
}\mathrel{\mathop:}=\mathcal{T}^{\delta }\left( \varphi \right)$,
as we know that $y^{\delta }$ converges to $\mathcal{T}\left(
\varphi \right)$ strongly in $H_{0}^{1}\left( \Omega \right)$, we
pass to the limit in (\ref{born_y_delta}).
\end{proof}

\begin{theo}
\label{theotocon} $\mathcal{T}$ is continuous from $\mathcal{U}$
endowed with the sequential weak topology of $H^{2}\left( \Omega
\right)$ to $H_{0}^{1}\left( \Omega \right) $ endowed with the
sequential weak topology.
\end{theo}

\begin{proof}
Let $\varphi _{k}$ be a sequence that converges to $\varphi$
weakly in  $H^{2}\left( \Omega \right)$. Then the sequence
$\varphi _{k}$ converges strongly in $H_{0}^{1}\left( \Omega
\right)$. We set $y_{k}\mathrel{\mathop:}=\mathcal{T}\left(
\varphi _{k}\right)$. Let $v$ in $\mathcal{K}\left( \varphi
\right)$ and set $v_{k}=\sup \left( v,\varphi _{k}\right)$ $\in$
$\mathcal{K}\left( \varphi _{k}\right)$. The sequence $v_{k}$
converges
to $v$ strongly in $H_{0}^{1}\left( \Omega \right)$. As $y_{k}\mathrel{\mathop:}=\mathcal{T}%
\left( \varphi _{k}\right) $, we get  $a\left(
y_{k},v_{k}-y_{k}\right) \geq \left( f,v_{k}-y_{k}\right)$, i.e.

\begin{equation*}
a\left( y_{k},y_{k}\right) \leq a\left( y_{k},v_{k}\right) -\left(
f,v_{k}-y_{k}\right).
\end{equation*}

By Corollary \ref{cortocon} the sequence $y_{k}$ is bounded and
weakly converges in $ H_{0}^{1}\left( \Omega \right)$ to some $y$.
Using the lower semi-continuity of $a$, the previous relation
gives

\begin{equation*}
a\left( y,y\right) \leq a\left( y,v\right) -\left( f,v-y\right).
\end{equation*}

As $y_{k}\geq \varphi _{k}$, this implies that  $y\geq \varphi$,
therefore $y\mathrel{\mathop:}=\mathcal{T}\left( \varphi \right)$.
\end{proof}

We obtain the main result of this section.
\begin{theo}
Problem $\left( \mathcal{P}\right)$ admits at least one optimal
solution $\varphi$ $\in$ $\U$.
\end{theo}

\begin{proof}
Let $\varphi _{k}$ a minimizing sequence. As $J\left( \varphi
_{k}\right) $ is bounded, $\varphi _{k}$ is $H^{2}$-bounded and
converges to $\varphi $ weakly in $H^{2}\left( \Omega \right) $.
By Theorem \ref{theotocon}, the sequence
$y_{k}\mathrel{\mathop:}=\mathcal{T}\left( \varphi _{k}\right)$
converges to $y^{\ast}\mathrel{\mathop:}=\mathcal{T}\left( \varphi
^{\ast }\right)$ weakly in $H_{0}^{1}\left( \Omega \right)$ and
using the norms semi-continuity we obtain

\begin{equation*}
J\left( \varphi ^{\ast
}\right)\mathrel{\mathop:}=\tfrac{1}{2}\int\nolimits_{\Omega
}\left( \mathcal{T}\left( \varphi ^{\ast }\right) -z\right)
^{2}d{x}+\tfrac{\nu}{2}\left( \int\nolimits_{\Omega }\left( \Delta
\varphi ^{\ast }\right) ^{2}d{x}\right) \leq
\underset{k\rightarrow \infty }{\lim \inf }J\left( \varphi
_{k}\right) =\inf \left( \mathcal{P}\right).
\end{equation*}
\end{proof}

\setcounter{equation}{0}

%-----------------------------------------------------------------------------------------

\subsection{An Approximated problem $\left( \mathcal{P}_{\delta}\right)
$}

%------------------------------------------------------------------------------------------

We use a trick of Barbu \cite{Bar1}, and add adapted penalization
terms to the approximated functional $J_{\delta}$ (here we add
$\tfrac{1}{2} \| \varphi-\varphi^{\ast} \|_{0}^{2}$) to force the
relaxed obstacle family $\varphi$ to converge to a desired
solution $\varphi^{\ast}$ of $(\mathcal{P})$ . So for any $\delta
>0$, we define

\begin{equation*}
J_{\delta }\left( \varphi \right)
\mathrel{\mathop:}=\tfrac{1}{2}\left[ \int\nolimits_{\Omega
}\left( \mathcal{T}^{\delta }\left( \varphi \right) -z\right)
^{2}d{x}+ \nu\left( \int\nolimits_{\Omega }\left( \Delta \varphi
\right) ^{2}d{x}\right) +\left\Vert \varphi -\varphi ^{\ast
}\right\Vert_{\Ld}^{2}\right].
\end{equation*}

The approximated optimal control problem  $\left(
\mathcal{P}^{\delta}\right)$ stands
\begin{equation}
\min \left\{ J_{\delta }\left( \varphi \right) ,\ \varphi \in
\mathcal{U}_{ad}\right\} \tag{${\mathcal{P}}^{\delta }$}.
\end{equation}

\begin{theo}
\label{conv} Problem $\left( \mathcal{P}^{\delta }\right)$ admits
at least one solution $\varphi ^{\delta }$. Moreover, when
$\delta$ go to $0$, the family $\varphi^{\delta}$ converges to
$\varphi^{\ast}$ weakly in $H^{2}\left( \Omega \right)$, and
$y^{\delta }\mathrel{\mathop:}=\mathcal{T}^{\delta }\left( \varphi
^{\delta }\right)$ converges to $y^{\ast
}\mathrel{\mathop:}=\mathcal{T}\left( \varphi ^{\ast }\right)$,
strongly in $H_{0}^{1}\left( \Omega \right)$.
\end{theo}

\begin{proof}
The functional $J_{\delta }$ is coercive, and lower
semi-continuous on $\mathcal{U}$. Therefore, the problem $\left(
\mathcal{P}^{\delta }\right) $ admits at least one solution
$\varphi ^{\delta }$. We set $y^{\delta
}\mathrel{\mathop:}=\mathcal{T}^{\delta }\left( \varphi ^{\delta
}\right)$, and note, that for any $\delta
>0$,

\begin{equation}
J_{\delta }\left( \varphi ^{\delta }\right) \leq J_{\delta }\left(
\varphi
^{\ast }\right)\mathrel{\mathop:}=\tfrac{1}{2}\left[ \int\nolimits_{\Omega }\left( \mathcal{T}%
^{\delta }\left( \varphi ^{\ast }\right) -z\right) ^{2}d{ x}+
\nu\left(
\int\nolimits_{\Omega }\left( \Delta \varphi ^{\ast }\right) ^{2}d{x}\right) %
\right].  \label{ineg_J}
\end{equation}

By Theorem \ref{ydcony}, we know that $\mathcal{T}^{\delta }\left(
\varphi ^{\ast }\right)$ converges to $\mathcal{T}\left( \varphi
^{\ast }\right)$ strongly in $H_{0}^{1}\left( \Omega \right)$, so
that $J_{\delta }\left( \varphi ^{\ast }\right)$ converges to
$J\left( \varphi ^{\ast }\right)$ as $\delta \rightarrow 0$.
Consequently, there exist $\delta _{0}>0$ and a constant $j^{\ast
}$, such that

\begin{equation*}
\forall ~\delta \leq \delta _{0},\quad J_{\delta }\left( \varphi
^{\delta }\right) \leq j^{\ast }<+\infty .
\end{equation*}

Consequently $\varphi ^{\delta }$ is $H^{2}$-bounded uniformly,
for any $\delta \leq \delta _{0}$. We use the Theorem
\ref{ydcony}, we get $\varphi^{\delta}$ converge to
$\widetilde{\varphi}$ weakly in $H^{2}\left(\Omega\right)$ and
strongly in $\Huz$ and $y^{\delta}$ converge to
$\widetilde{y}\mathrel{\mathop:} =
\mathcal{T}\left(\widetilde{\varphi}\right)$ strongly in $\Huz$.
As $\mathcal{U}_{ad}$ is weakly closed, we have
$\widetilde{\varphi}$ in $\mathcal{U}_{ad}$. By \eqref{ineg_J} and
the lower semi-continuity of $J_{\delta}$, we get

\begin{equation*}
\begin{split}
J\left( \widetilde{\varphi }\right)+ \tfrac{1}{2}\left\Vert
\widetilde{\varphi }-\varphi ^{\ast }\right\Vert_{0}^{2} & \leq
\underset{\delta \rightarrow 0}{~\lim \inf}J_{\delta }\left(
\varphi ^{\delta }\right) \\
& \leq  \underset{\delta \rightarrow 0}{~\lim \sup}J_{\delta
}\left(
\varphi ^{\delta }\right)  \\
& \leq  \underset{\delta \rightarrow 0}{~\lim }J_{\delta }\left(
\varphi ^{\ast }\right) \leq \underset{\delta \rightarrow 0}{~\lim
}J\left( \varphi^{\ast }\right)  \\
& \leq  J\left( \widetilde{\varphi }\right) .
\end{split}
\end{equation*}
This yields that $\left\Vert \widetilde{\varphi }-\varphi ^{\ast
}\right\Vert_{0}^{2} \leq 0$, then $\widetilde{\varphi} =
\varphi^{\ast}$ and $ \underset{\delta \rightarrow 0}{\lim
}J_{\delta }\left( \varphi ^{\delta }\right) =J\left( \varphi
^{\ast }\right).$

In addition this proves that any clusters points of
$\varphi^{\delta}$ is equal to $\varphi^{\ast}$, so that the whole
family converges.

\end{proof}

\setcounter{equation}{0}

%----------------------------------------------------------------------

\subsection{Optimality conditions for problem $\left( \mathcal{P}^{\protect\delta }\right)$}

%-----------------------------------------------------------------------

We give first necessary optimality conditions for problem $\left(
\mathcal{P}^{\delta }\right)$. Les us recall the following result
on the Gâteaux-derivative of the operator $\mathcal{T}^{\delta }$
[1].

\begin{lem}
The mapping  $\mathcal{T}^{\delta }$ is Gâteaux-derivative at any
$\varphi $ in $\mathcal{U}_{ad}$:

\begin{equation*}
\forall \xi \in H_{0}^{1}\left( \Omega \right) ,\quad
\frac{\mathcal{T}^{\delta }\left( \varphi +\tau \xi \right)
-\mathcal{T}^{\delta }\left( \varphi \right)
}{\tau}\overset{w}{\rightharpoonup }v^{\delta
},~\text{in}~H_{0}^{1}\left( \Omega \right),~\text{when}~\tau
\rightarrow 0,
\end{equation*}
where $v^{\delta }$ is the solution of the following equation
\begin{equation*}
Av^{\delta }+ \beta_{\delta } ^{\prime }\left( y^{\delta }-\varphi
\right)v^{\delta }=\beta_{\delta } ^{\prime }\left( y^{\delta
}-\varphi \right)\xi \quad \mathrm{in}~\Omega ,\quad v^{\delta
}=0\quad \mathrm{on}~\partial \Omega.
\end{equation*}

\end{lem}
\begin{proof}
See \cite{Adam}.
\end{proof}

We define the approximate adjoint state $p^{\delta }$ in
$H_{0}^{1}\left( \Omega \right) $ as the solution of the following
adjoint equation

\begin{equation*}
\label{eq_adjon1} A^{\ast}p^{\delta }+\beta_{\delta } ^{\prime
}\left( y^{\delta }-\varphi \right)p^{\delta }=y^{\delta }-z \quad
\mbox{in}~\Omega,\quad p^{\delta } = 0\quad \mbox{on}~\partial
\Omega,
\end{equation*}

where $A^{\ast}$ is the adjoint operator of $A$. As $\varphi
^{\delta }$ is the solution  of the problem $\left(
\mathcal{P}^{\delta }\right) $, we get
\begin{equation*}
\forall \varphi \in \mathcal{U}_{ad},\quad \frac{d}{dt}J_{\delta
}\left( \varphi ^{\delta }+t\left( \varphi -\varphi ^{\delta
}\right) \right) _{\mid t=0}\geq 0.
\end{equation*}

That is

\begin{center}
\begin{equation*}
\forall \varphi \in \mathcal{U}_{ad},\quad\int\nolimits_{\Omega
}\left( \chi ^{\delta }\left( y^{\delta }-z\right) +\nu \Delta
\varphi ^{\delta }\Delta \left( \varphi -\varphi ^{\delta }\right)
\right) ~d{x}+\int\nolimits_{\Omega }\left( \varphi ^{\delta
}-\varphi ^{\ast }\right) \left( \varphi -\varphi ^{\delta
}\right)~d{x}\geq 0,
\end{equation*}
\end{center}

where $\chi ^{\delta } \in H_{0}^{1}\left( \Omega \right)$ and
satisfies

\begin{equation*}
A\chi ^{\delta }+\beta_{\delta } ^{\prime }\left( y^{\delta
}-\varphi ^{\delta }\right)\chi ^{\delta }=\beta_{\delta} ^{\prime
}\left( y^{\delta }-\varphi ^{\delta }\right) \left( \varphi
-\varphi ^{\delta }\right)\quad \text{in}~\Omega.
\end{equation*}

from the definition of $p^{\delta }$, we obtain

\begin{equation*}
\int\nolimits_{\Omega }\chi ^{\delta
}A^{\ast}p^{\delta}d{x}+\int\nolimits_{\Omega }\beta_{\delta }
^{\prime }\left( y^{\delta }-\varphi ^{\delta }\right) p^{\delta
}\chi ^{\delta }d{x} +\nu\int\nolimits_{\Omega }\Delta \varphi
^{\delta}\Delta \left( \varphi -\varphi ^{\delta }\right)\,d{x}
+\int\nolimits_{\Omega }\left( \varphi ^{\delta }-\varphi ^{\ast
}\right) \left( \varphi -\varphi ^{\delta }\right)\,d{x}\geq 0,
\end{equation*}

where $A^{\ast }$ denotes the adjoint operator of de $A$. Then

\begin{equation*}
\int\nolimits_{\Omega }A\chi ^{\delta
}p^{\delta}d{x}+\int\nolimits_{\Omega }\beta_{\delta} ^{\prime
}\left( y^{\delta }-\varphi ^{\delta }\right) p^{\delta }\chi
^{\delta }\,d{x} +\nu \int\nolimits_{\Omega }\Delta \varphi
^{\delta }\Delta \left( \varphi -\varphi ^{\delta
}\right)\,d{x}+\int\nolimits_{\Omega }\left( \varphi ^{\delta
}-\varphi ^{\ast }\right) \left( \varphi -\varphi ^{\delta
}\right)\, d{x}\geq 0,
\end{equation*}

we obtain

\begin{equation*}
\int\nolimits_{\Omega }\beta_{\delta }^{\prime }\left( y^{\delta
}-\varphi ^{\delta }\right) p^{\delta }\left( \varphi -\varphi
^{\delta }\right)\,d{x}+\nu \int\nolimits_{\Omega }\Delta \varphi
^{\delta }\Delta \left( \varphi -\varphi ^{\delta }\right)\,
d{x}+\int\nolimits_{\Omega }\left( \varphi ^{\delta }-\varphi
^{\ast }\right) \left( \varphi -\varphi ^{\delta }\right)\,
d{x}\geq 0.
\end{equation*}

In the sequel, we set

\begin{equation}
\label{mu} \mu ^{\delta }\mathrel{\mathop:}=\beta_{\delta
}^{\prime }\left( y^{\delta }-\varphi ^{\delta }\right)p^{\delta }
\in L^{2}\left( \Omega \right).
\end{equation}

Finally, we obtain

\begin{theo}
\label{theo_opti_sys_Ld}
If $\varphi ^{\delta }$ is an optimal solution of $\left( \mathcal{P}^{\delta }\right) $ and $y^{\delta }\mathrel{\mathop:}=\mathcal{T}%
^{\delta }\left( \varphi ^{\delta }\right)$, there exists
$p^{\delta }$ in $H^{2}(\Omega)\cap H_{0}^{1}(\Omega)$ and $\mu
^{\delta }$ in $L^{2}\left( \Omega \right) $ such that the
following system holds

\begin{subequations}
\label{opti_sys_gener}
\begin{equation}  \label{equastat}
Ay^{\delta }+\beta_{\delta }\left( y^{\delta }-\varphi ^{\delta
}\right)=f \quad \text{in}~\Omega,\quad y^{\delta }=0\quad
\text{on}~\partial \Omega ,
\end{equation}
\begin{equation}  \label{equaadjoint}
A^{\ast }p^{\delta}+\mu ^{\delta }=y^{\delta }-z \quad
\text{in}~\Omega ,\quad p^{\delta }=0\quad \text{on}~\partial
\Omega ,
\end{equation}
\begin{equation}
\left( \mu ^{\delta }+\varphi ^{\delta }-\varphi ^{\ast },\varphi
-\varphi ^{\delta }\right)+\nu \left( \Delta \varphi ^{\delta
},\Delta \left( \varphi -\varphi ^{\delta }\right) \right)\geq
0,\quad \forall \varphi \in \mathcal{U}_{ad}. \label{ineq_proj}
\end{equation}
\end{subequations}
\end{theo}

In the case
$\mathcal{U}_{ad}\mathrel{\mathop:}=L^{2}\left(\Omega\right)$,  we
make this optimality system more precise.

Let $\chi $ in $\mathcal{U}$ and choose  $\varphi =\varphi
^{\delta }\pm \chi$; by the equation (\ref{ineq_proj}), we obtain

\begin{equation}
\left( \mu ^{\delta }+\varphi ^{\delta }-\varphi ^{\ast },\chi
\right)+\nu \left( \Delta \varphi ^{\delta },\Delta \chi
\right)=0, \quad \forall \chi \in \mathcal{U}.  \label{egal}
\end{equation}

Set $h^{\delta }=\Delta \varphi ^{\delta }$ in $L^{2}\left( \Omega
\right)$, so that for any $\chi $ in $\mathscr{D}\left( \Omega
\right) $, the relation (\ref{egal}) gives

\begin{equation*}
\left( \mu ^{\delta }+\varphi ^{\delta }-\varphi ^{\ast },\chi
\right)+\nu \left( h^{\delta },\Delta \chi \right)=0~\mbox{(in the
distribution sens)},
\end{equation*}

that is

\begin{equation*}
-\nu \Delta h^{\delta }=\mu ^{\delta }+\varphi ^{\delta }-\varphi ^{\ast }%
\in \mathscr{D}^{\prime }\left( \Omega \right).
\end{equation*}

Using the same techniques as in \cite{Bergou1}, we deduce that
$h_{\mid \partial \Omega }^{\delta }=0$. Consequently, $h^{\delta
}$ $\in$ $\mathcal{U}$, and it is the unique solution of

\begin{equation*}
-\nu \Delta h^{\delta }=\mu ^{\delta }+\varphi ^{\delta }-\varphi
^{\ast }\text{ in } L^{2}\left( \Omega \right) ,\quad h^{\delta
}=0\quad \text{on}~\partial \Omega.
\end{equation*}

The last relation may be written as

\begin{equation*}
-\nu \left(\Delta ^{2}\varphi ^{\delta },u \right)=\left(\mu
^{\delta },u\right)
-\left(\varphi^{\delta}-\varphi^{\ast},u\right)\quad \text{in}~
\Omega ,\quad \varphi ^{\delta }=0\quad \text{on }\partial \Omega.
\label{projequa1}
\end{equation*}

Thanks to Green's formula, the previous relation reads

\begin{multline*}
-\nu \int\nolimits_{\Omega }\Delta ^{2}\varphi\,u\,d{x}- \left(\mu
^{\delta
},u\right)-\left(\varphi^{\delta}-\varphi^{\ast},u\right)=\\=\int\nolimits_{\Omega
}\Delta
\varphi \Delta u\, d{x}-\int\nolimits_{\Gamma }\left( \Delta \varphi \frac{\partial u}{%
\partial \eta }-u\frac{\partial \Delta \varphi}{\partial \eta
}\right)-\left(\mu ^{\delta
},u\right)-\left(\varphi^{\delta}-\varphi^{\ast},u\right).
\end{multline*}

So $\Delta \varphi$ vanishes on the boundary $\partial \Omega$,
and we conclude that $\varphi^{\delta}$ belongs to
$\mathcal{W}\mathrel{\mathop:}= \{ u\,|\,u \in H^{2}(\Omega) \cap
H^{1}_{0}(\Omega)\, \text{et}\, \Delta u_{|\partial \Omega} = 0
\}$. Finally we have:

\begin{corol}
Assume conditions of Theorem \ref{theo_opti_sys_Ld} are fulfilled,
and $\mathcal{U}_{ad} = \mathcal{U}$, then the optimality system
$({\mathcal{S}}^{\delta })$ reads
\begin{subequations}
\label{optimalsysvar}
\begin{equation}  \label{statequavar}
a\left( y^{\delta },v\right) +\left( \beta _{\delta }\left(
y^{\delta }-\varphi ^{\delta }\right) ,v\right)=\left(
f,v\right),\quad \forall v\in \mathcal{U},
\end{equation}
\begin{equation}  \label{adjstatequavar}
a^{\ast}\left( p^{\delta },w\right) +\left( \mu ^{\delta
},w\right)=\left( y^{\delta }-z,w\right),\quad \forall w\in
\mathcal{U},
\end{equation}
\begin{equation}  \label{projequavar}
\nu \left( \Delta ^2\varphi ^{\delta },u \right) - \left(
 \mu ^{\delta }, u\right) = \left(\varphi^{\delta}-\varphi^{\ast}, u\right) ,\quad \forall u \in \mathcal{W}.
\end{equation}
\end{subequations}

Here $a^{\ast}$ denotes the adjoint form of $a$ (associated with
the adjoint operator $A^{\ast}$).

\end{corol}

\setcounter{equation}{0}
%------------------------------------------------------------------------------------------------
\section{First order necessary optimality conditions for $\left(\mathcal{P}\right)$}
%------------------------------------------------------------------------------------------------

In this section, we have to estimate $p^{\delta}$, and gives more
convergence results, then we may pass to the limit in the system
\eqref{optimalsysvar} as $\delta \rightarrow 0$.

\begin{theo}
When $\delta $ goes to  $0$, $p^{\delta }$ converges to $p^{\ast
}$ weakly in $H_{0}^{1}\left( \Omega \right)$ and $\mu ^{\delta }$
converges to $\mu^*$ weakly satr in $H^{-1}\left( \Omega \right)
\cap \mathcal{M}\left(\Omega\right)$, and

\begin{equation*}
\left\langle \mu ^{\ast},p^{\ast }\right\rangle \geq 0.
\end{equation*}

where $\mathcal{M}\left(\Omega\right)$ is the set of all regular
signed measures in $\Omega$.

\end{theo}

\begin{proof}
Using \eqref{adjstatequavar}, we obtain
\begin{equation}
a^{\ast }\left( p^{\delta },p^{\delta }\right)
+\int\nolimits_{\Omega }{\beta_{\delta }^{\prime }\left( y^{\delta
}-\varphi ^{\delta }\right) }\left( p^{\delta }\right)
^{2}dx=\left( y^{\delta }-z,p^{\delta }\right).
\label{equ_adjoin_var}
\end{equation}

As $\beta ^{\prime }\left(\cdot \right)\geq 0$, and thanks to
hypothesis \eqref{conbilfora} and \eqref{coebilfora}, we get
\begin{equation*}
\left\Vert p^{\delta }\right\Vert_{\Hu} \leq C\left\Vert y^{\delta
}-z\right\Vert_{\Hu},
\end{equation*}

which implies that $p^{\delta }$ converges to $p^{\ast }$ weakly
in $H_{0}^{1}\left( \Omega \right)$. Consequently  $A^{\ast}
p^{\delta }$ is uniformly bounded  in $H^{-1}\left( \Omega \right)
$ and
\begin{equation}
\label{pdconp} \mu ^{\delta }=-A^{\ast}p^{\delta }+y^{\delta }-z.
\end{equation}

Let $\gamma _{\varepsilon }\in \mathscr{C}^{1}\left(
\mathbb{R}\right) $ be a family of smooth approximations to the
sign function and satisfy the following \cite{Yuquan}:

\begin{equation*}
\gamma _{\varepsilon }^{\prime }\left( r\right) \geq 0\quad
\forall r\in \mathbb{R},
\end{equation*}%

and

\begin{equation*}
\gamma _{\varepsilon }^{\prime }\left( r\right)
\mathrel{\mathop:}=
\begin{cases}
1 & \mathrm{if\ }  r>\varepsilon \\
0 & \mathrm{if\ }  r=0 \\
-1 & \mathrm{if\ }  r<-\varepsilon,%
\end{cases}%
\end{equation*}%

Then we can multiply \eqref{pdconp} by $\gamma _{\varepsilon
}\left( p^{\delta }\right) $ and integrate it over $\Omega $. As a
result, we get

\begin{equation*}
\int\nolimits_{\Omega }\mu ^{\delta }\gamma _{\varepsilon }\left(
p^{\delta }\right) dx\leq C.
\end{equation*}

Letting $\varepsilon \rightarrow 0$, we have

\begin{equation*}
\left\Vert \mu ^{\delta }\right\Vert _{L^{1}\left( \Omega \right)
}\leq C.
\end{equation*}

Hence $\mu ^{\delta }$ is bounded in $L^{1}\left( \Omega \right) $
and consequently it is also bounded in $\mathcal{M}\left( \Omega
\right) $, thus, $\mu ^{\delta }$ converge to $\mu^{\ast}$ weakly
star in $H^{-1}\left( \Omega \right) \cap \mathcal{M}\left( \Omega
\right)$, such that
\begin{equation}
A^{\ast} p^{\ast }+\mu ^{\ast }=y^{\ast }-z\quad \text{in}\quad
\Omega ,\quad p^{\ast }=0\quad \text{on}~\partial \Omega.
\label{equa_adjoin}
\end{equation}
As $0\leq \beta ^{\prime }\leq 1$, by using (\ref{equ_adjoin_var})
we get
\begin{equation*}
a^{\ast }\left( p^{\delta },p^{\delta }\right) \leq \left(
y^{\delta }-z,p^{\delta }\right).
\end{equation*}

And by the lower semi-continuity of $a^{\ast }$

\begin{equation*}
\langle A^{\ast}p^{\ast}, p^{\ast}\rangle = a^{\ast }\left(
p^{\ast },p^{\ast }\right) \leq \underset{\delta \rightarrow
\infty }{\lim \inf }\left( y^{\delta }-z,p^{\delta }\right)=\left(
y^{\ast }-z,p^{\ast }\right).
\end{equation*}

From (\ref{equa_adjoin}), we obtain
\begin{equation*}
0\leq \left\langle A^{\ast} p^{\ast },p^{\ast }\right\rangle
=\left( y^{\ast }-z,p^{\ast }\right)-\left\langle \mu ^{\ast
},p^{\ast }\right\rangle \leq \left( y^{\ast }-z,p^{\ast }\right),
\end{equation*}

so that

\begin{equation*}
\left\langle \mu ^{\ast},p^{\ast }\right\rangle \geq 0.
\end{equation*}

\end{proof}

In the sequel, we set $\xi^{\delta } \mathrel{\mathop:} = \beta
_{\delta }\left( y^{\delta }-\varphi ^{\delta }\right)$, then we
obtain the following results.

\begin{theo}
\label{theo_conv_xi} When  $\delta $ goes to $0$, $\xi^{\delta }$
converges to $\xi^{\ast}$ weakly in  $\Ld$, where $\xi^{\ast}$ is
negative and the state equation \eqref{equastat} gives

\begin{equation*}%
Ay^{\ast}+\xi^{\ast}=f.
\end{equation*}
\end{theo}

\begin{proof}
From \eqref{statequavar}, we obtain with
$v=\beta_{\delta}\left(y^{\delta}-\varphi^{\delta}\right)$

\begin{equation*}
\label{eq_eta_exp}%
\left\langle A\left( y^{\delta }-\varphi ^{\delta }\right) ,\beta
_{\delta }\left( y^{\delta }-\varphi ^{\delta }\right)
\right\rangle +\left( \beta _{\delta }\left( y^{\delta }-\varphi
^{\delta }\right) ,\beta _{\delta }\left( y^{\delta }-\varphi
^{\delta }\right) \right) =\left( f-A\varphi ^{\delta },\beta
_{\delta }\left( y^{\delta }-\varphi ^{\delta }\right) \right).
\end{equation*}

For the seek of simplicity, we set $r^{\delta
}\mathrel{\mathop:}=y^{\delta}-\varphi ^{\delta }$; with
\eqref{bilfora}, this gives

\begin{multline}
\label{stat_eq_expli} \sum\limits_{i,j=1}^{n}\int\nolimits_{\Omega
}a_{ij} \frac{\partial r^{\delta } }{\partial {x}_{i}}
\frac{\partial r^{\delta } }{\partial {x}_{j}
}\beta_{\delta}^{\prime }\left(r^{\delta}\right) d{
x}+\int\nolimits_{\Omega }a_{0} r^{\delta } \beta_{\delta} \left(
r^{\delta }\right) d{x}+\left\Vert \beta _{\delta }\left(
r^{\delta }\right)
\right\Vert^{2}_{\Ld}+\\+\sum\limits_{i=1}^{n}\int\nolimits_{\Omega
}a_{i} \frac{\partial r^{\delta } }{\partial
{x}_{i}}\beta_{\delta} \left(r^{\delta}\right) d{x} =\left(
f-A\varphi ^{\delta },\beta _{\delta}\left(r^{\delta}\right)
\right).
\end{multline}

With hypothesis $(\mathbf{H})$, we get

\begin{equation}
\begin{split}
\label{ellip_cond_expli}
\sum\limits_{i,j=1}^{n}\int\nolimits_{\Omega }a_{ij} \frac{%
\partial r^{\delta} }{\partial {x}_{i}}%
\frac{\partial r^{\delta} }{\partial {x}_{j}%
}\beta _{\delta }^{\prime }\left( r^{\delta}\right)
d{x}&\geq \int\nolimits_{\Omega }m  \sum\limits_{i=0}^{n}\left( \frac{%
\partial r^{\delta} }{\partial {x}_{i}}\right)^{2}\beta _{\delta }^{\prime }\left( r^{\delta}\right)
d{x}.\\ & \geq 0.
\end{split}
\end{equation}

From \eqref{stat_eq_expli}, \eqref{ellip_cond_expli} and
$(\mathbf{H})$, we obtain

\begin{equation*}
\begin{split}
\left\Vert \beta _{\delta }\left( r^{\delta}\right) \right\Vert
_{\Ld}^{2}& \leq \left\Vert f-A\varphi^{\delta} \right\Vert
_{\Ld}\left\Vert \beta _{\delta }\left( r^{\delta}\right)
\right\Vert _{\Ld}+\sum\limits_{i=1}^{n}\left\Vert a_{i}
\right\Vert _{L^{\infty }\left( \Omega \right) }\left( \left\Vert
\nabla r^{\delta} \right\Vert _{L^{2}\left( \Omega \right)
}\left\Vert \beta _{\delta }\left(r^{\delta}\right)
\right\Vert _{\Ld}\right) \\
& \leq \max \left\{ 1,\sum\limits_{i=1}^{n}\left\Vert a_{i}
\right\Vert _{L^{\infty }\left( \Omega \right) }\right\} \left(
\left\Vert f-A\varphi^{\delta} \right\Vert _{\Ld}+\left\Vert
\nabla r^{\delta} \right\Vert _{\Ld}\right) \left\Vert \beta
_{\delta }\left( r^{\delta}\right) \right\Vert _{\Ld}.
\end{split}
\end{equation*}

Finally, we get

\begin{equation*}
\left\Vert \beta _{\delta }\left( y^{\delta }-\varphi ^{\delta
}\right) \right\Vert _{\Ld}\leq \alpha \left( \left\Vert
f-A\varphi^{\delta} \right\Vert _{\Ld}+\left\Vert \nabla \left(
y^{\delta }-\varphi ^{\delta }\right) \right\Vert _{\Ld}\right),
\end{equation*}

where $\alpha \mathrel{\mathop:}= \max \left\{
1,\sum\limits_{i=1}^{n}\left\Vert a_{i} \right\Vert _{L^{\infty
}\left( \Omega \right) }\right\}$, so that

\begin{equation*}
\label{est_beta} \left\Vert \beta _{\delta }\left( y^{\delta
}-\varphi ^{\delta }\right) \right\Vert _{\Ld}\leq \alpha \left(
\left\Vert f-A\varphi^{\delta} \right\Vert _{\Ld}+\left\Vert
y^{\delta }-\varphi ^{\delta }\right\Vert_{\Hu}\right).
\end{equation*}

Since $\varphi ^{\delta }$ and $y^{\delta }$ are respectively
bounded in $H^{2}\left( \Omega \right) \cap H_{0}^{1}\left( \Omega
\right)$ and  $\Huz$, we deduce that $\xi^{\delta}$ is bounded in
$\Ld$, by passing to the limit where $\delta \rightarrow 0$, we
obtain that $\xi ^{\delta}$ converge to $\xi^{\ast}$ weakly in
$\Ld$. Passing to the limit in \eqref{equastat}, gives

\begin{equation*}
 A y^{\ast}+\xi^{\ast}=f.
\end{equation*}

where $\xi^{\ast}$, is negative and we get $y^{\ast} \in \Hd \cap
\Huz$.
\end{proof}

\begin{corol}

As $\varphi^{\ast}$ is in $\Hd \cap \Huz$,
$y^{\ast}\mathrel{\mathop:}=\mathcal{T}^{\ast}\left(
\varphi^{\ast} \right)$ belongs to $\Hd \cap \Huz$.
\end{corol}

\begin{proof}
As $\xi^{\ast}$
 and $f$ belongs to $\Ld$, then $Ay^{\ast} \in \Ld$ and $y^{\ast}\in
 \Hd$.
\end{proof}

Now, we give some Lemmas, the proof the below Theorem
\ref{theo_syst_conv}.

\begin{lem}
\label{lem_xi_p_0} When $\delta$ goes to $0$, $\left( \mu^{\delta
},\left( y^{\delta }-\varphi ^{\delta }\right) ^{+}\right)
\rightarrow \left\langle \mu^{\ast} ,y^{\ast}-\varphi^{\ast}
\right\rangle$, and $\left\langle \mu^{\ast}
,y^{\ast}-\varphi^{\ast} \right\rangle =0$.
\end{lem}

\begin{proof}
By the definition of $\beta$ and $\mu^{\delta}$ \eqref{mu}, we get

\begin{equation*}
\left( \mu^{\delta },\left( y^{\delta }-\varphi ^{\delta }\right)
^{+}\right) =0,
\end{equation*}

where $v^{+}\mathrel{\mathop:}=\max \left\{ 0,v\right\}$, by
Theorem \ref{ydcony} $\left( y^{\delta }-\varphi ^{\delta
}\right)^{+}$ converges strongly to $\left(y^{\ast}-\varphi^{\ast}
\right)$ in $H_{0}^{1}\left( \Omega \right) $,

then

\begin{equation*}
\left( \mu^{\delta },\left( y^{\delta }-\varphi ^{\delta }\right)
^{+}\right) \rightarrow \left\langle \mu^{\ast}
,y^{\ast}-\varphi^{\ast} \right\rangle,
\end{equation*}

and

\begin{equation*}
\left\langle \mu^{\ast} ,y^{\ast}-\varphi^{\ast} \right\rangle =0.
\end{equation*}

\end{proof}

\begin{lem}
\label{lem_xi_p_0} When $\delta$ goes to $0$, $\left( \xi ^{\delta
},p^{\delta }\right) \rightarrow \left\langle \xi
^{\ast},p^{\ast}\right\rangle$, and $\left\langle \xi
^{\ast},p^{\ast}\right\rangle=0$.
\end{lem}

\begin{proof}
As befor, we set $r^{\delta }\mathrel{\mathop:}=y^{\delta}-\varphi
^{\delta }$, so that
$\xi^{\delta}=\beta_{\delta}\left(r^{\delta}\right)$. From the
definition of $\beta$ and $\beta ^{\prime}$, we get respectively

\begin{equation*}
\begin{split}
\left( \xi^{\delta },p^{\delta }\right)=
\left(\beta_{\delta}\left(r^{\delta}\right),p^{\delta}\right)&=
\\&=\tfrac{1}{\delta
}\left[ \int\nolimits_{\left\{ r^{\delta}\leq
-\tfrac{1}{2}\right\} }\left( r^{\delta
} +\tfrac{1}{4}\right) p^{\delta }d{x}%
-\int\nolimits_{\left\{ -\tfrac{1}{2}\leq r^{\delta }\leq
0\right\} }\left( r^{\delta
}\right) ^{2}p^{\delta }d{%
x}\right],
\end{split}
\end{equation*}

and

\begin{equation*}
\begin{split}
\left( \mu^{\delta
},r^{\delta}\right)=\left(\beta_{\delta}^{\prime}\left(r^{\delta}p^{\delta}\right),p^{\delta}\right)
&=\\&=\tfrac{1}{\delta }\left[ \int\nolimits_{\left\{
r^{\delta}\leq -\tfrac{1}{2}\right\} }p^{\delta } r^{\delta}
d{x}-2\int\nolimits_{\left\{ -\tfrac{1}{2}\leq r^{\delta}\leq
0\right\} }\left( r^{\delta}\right) ^{2}p^{\delta }d{x}\right].
\end{split}%
\end{equation*}

from that, we get

\begin{equation*}
\left( \xi ^{\delta },p^{\delta }\right) -\tfrac{1}{2}%
\left( \mu ^{\delta },r^{\delta}\right) =
\tfrac{1}{\delta }\int\nolimits_{\left\{ r^{\delta}\leq -%
\tfrac{1}{2}\right\} }\tfrac{1}{2}\left( r^{\delta}
+\tfrac{1}{2}\right) p^{\delta }d{\pmb x},
\end{equation*}

Then, we obtain

\begin{equation*}
\left\vert \left( \xi ^{\delta },p^{\delta }\right) -\tfrac{1}{2}%
\left( \mu ^{\delta },r^{\delta}\right) \right\vert \leq
\tfrac{1}{\delta }\left( \int\nolimits_{\left\{ r^{\delta}\leq
-\tfrac{1}{2}\right\} }\tfrac{1}{4}\left(
r^{\delta} +\tfrac{1}{2}\right)^{2}d{x}%
\right) ^{1/2}\left( \int\nolimits_{\left\{ r^{\delta}\leq
-\tfrac{1}{2}\right\} }\left( p^{\delta }\right) ^{2}d{x}\right)
^{1/2}.
\end{equation*}

As $\Hu \hookrightarrow L^{q}\left(\Omega\right)$ with $ 2 < q
\leq 6$, we have

\begin{equation*}
\begin{split}
\left( \int\nolimits_{\left\{ r^{\delta}\leq -\frac{1}{2}%
\right\}}(p^{\delta})^{2}d{x}\right) ^{1/2} & \leq
\left( \int\nolimits_{\left\{ r^{\delta}\leq -\frac{1}{2%
}\right\} }\left( p^{\delta }\right) ^{q}d{x}\right) ^{1/q}
\left( \int\nolimits_{\left\{ r^{\delta}\leq -\frac{1}{2%
}\right\}}d{x}\right) ^{(q-2)/2q}\\
&\leq C\left\Vert p^{\delta }\right\Vert _{\Hu} \left( \m \left\{
r^{\delta}\leq -\tfrac{1}{2}\right\} \right) ^{\left(q-2\right)/2q
},
\end{split}
\end{equation*}

where  ${\m} \left\{ \mathbb{A}\right\} $ is the measure of the
set $\mathbb{A}$. Then we write

\begin{equation}
\left\vert \left( \xi ^{\delta },p^{\delta }\right) -\tfrac{1}{2}
\left( \mu ^{\delta },r^{\delta}\right) \right\vert \leq
C\left\Vert p^{\delta }\right\Vert _{\Hu}\tfrac{1}{\delta }\left(
\int\nolimits_{\left\{ r^{\delta}\leq -\tfrac{1}{2}\right\}
}\tfrac{1}{4}\left(r^{\delta} +\frac{1}{2}\right)^{2}d{x}\right)
^{1/2}\left( \m \left\{ r^{\delta}\leq -\tfrac{1}{2}%
\right\} \right) ^{\left(q-2\right)/2q}. \label{est_cherch}
\end{equation}

We have

\begin{equation}
\tfrac{1}{\delta }\left( \int\nolimits_{\left\{ r^{\delta }\leq -\tfrac{1}{2}%
\right\} }\tfrac{1}{4}\left( r^{\delta }+\tfrac{1}{2}\right)
^{2}d{x}\right) ^{1/2}\leq
\tfrac{1}{\delta }\left( \int\nolimits_{\left\{ r^{\delta }\leq -\tfrac{1}{2}%
\right\} }\left( r^{\delta }+\tfrac{1}{4}\right) ^{2}d{x}\right)
^{1/2}. \label{est_r+05_lq_est_r+025}
\end{equation}

As  $\|\beta_{\delta} \left( r^{\delta }\right)\|_{\Ld}$ is
bounded, we deduce that

\begin{equation*}
\label{r+025_lq_M}
\tfrac{1}{\delta^{2}} \int\nolimits_{\left\{ r^{\delta }\leq -\tfrac{1}{2}%
\right\} }\left( r^{\delta }+\tfrac{1}{4}\right) ^{2} d{x}\leq C,
\end{equation*}

and

\begin{equation*}
\m \left\{ r^{\delta }\leq -\tfrac{1}{2}\right\} \leq C\delta
^{2}.
\end{equation*}%

then, by passing to the limit when $\delta$ goes to $0$, we have

\begin{equation*}
\underset{\delta \rightarrow 0}{\lim }\left(\m \left\{
r^{\delta}\leq -\tfrac{1}{2}\right\}\right)=0. \label{mes_r_lq_05}
\end{equation*}

Since $\|p^{\delta}\|_{\Hu}$ is bounded,  from \eqref{est_cherch}
and \eqref{est_r+05_lq_est_r+025}, we obtain

\begin{equation*}
\left\vert \left( \xi ^{\delta },p^{\delta }\right) -\tfrac{1}{2}%
\left( \mu ^{\delta },r^{\delta}\right) \right\vert \leq C\left(
\m \left\{r^{\delta}\leq -\tfrac{1}{2}\right\} \right)
^{\left(q-2\right)/2q},
\end{equation*}

and

\begin{equation*}
\underset{\delta \rightarrow 0}{\lim }\left\vert \left( \xi
^{\delta },p^{\delta }\right) -\tfrac{1}{2}\left( \mu ^{\delta
},r^{\delta}\right) \right\vert =0,
\end{equation*}

so when $\delta$ goes to $0$, $\left( \xi ^{\delta },p^{\delta
}\right)$ converges to $0$, and with Lemma \ref{lem_xi_p_0}, we
get

\begin{equation*}
\left\langle \xi^{\ast},p^{\ast}\right\rangle =0.
\end{equation*}
\end{proof}

Finally, from the previous convergence results, we obtain the main
result of this section

\begin{theo}
\label{theo_syst_conv} Let $\varphi ^{\ast }$, be an optimal
solution of problem $\left(\mathcal{P}\right)$. Then $\Delta
\varphi ^{\ast }$ belong to $H_{0}^{1}\left( \Omega \right) $ and
there exists $p^{\ast }$ in
 $H_{0}^{1}\left( \Omega \right)$, $\xi^{\ast} \leq 0$ in $\Ld$ and $\mu^{\ast }$  in $H^{-1}\left( \Omega \right)\cap \mathcal{M}\left(\Omega\right)$, such that the following optimality
 $(\mathcal{S})$ system holds

\begin{subequations}
\begin{equation}
A y^{\ast }+\xi ^{\ast }=f\quad \text{ in }\Omega ,\quad y^{\ast }=0\text{%
\quad on }\partial \Omega, \tag{1.a}
\end{equation}
\begin{equation}
A^{\ast} p^{\ast}+\mu ^{\ast }=y^{\ast }-z\quad \text{in }\Omega
,\quad p^{\ast }=0\text{\quad on }\partial \Omega, \tag{2.a}
\end{equation}
\begin{equation}
-\nu \Delta ^{2}\varphi ^{\ast }+\mu ^{\ast }=0\quad \text{in
}\Omega ,\quad \Delta \varphi ^{\ast }=\varphi ^{\ast
}=0\text{\quad on }\partial \Omega, \tag{3.a}
\end{equation}
\begin{equation}
\left\langle \mu ^{\ast },y^{\ast }-\varphi ^{\ast }\right\rangle
=0, \tag{4.a}
\end{equation}

\begin{equation}
\left\langle \xi ^{\ast },p^{\ast }\right\rangle=0, \tag{5.a}
\end{equation}

\begin{equation}
a^{*}\left(p^{\ast},p^{\ast}\right) -\left( z-y^{\ast },p^{\ast
}\right) \leq 0, \tag{6.a}
\end{equation}
\begin{equation}
\left\langle p^{\ast },\mu ^{\ast }\right\rangle \geq 0. \tag{7.a}
\end{equation}

\end{subequations}
\end{theo}

\begin{rem}
In \cite{Kunish}, Ito et Kunish had obtained the following
optimality condition system $(\widetilde{\mathcal{S}})$

\begin{subequations}

\begin{equation}
\tag{1.b} A y^{\ast}+\xi^{\ast} =f,\quad \xi^{\ast} =\max \left(
0,\xi^{\ast} +y^{\ast}-\varphi^{\ast} \right),
\end{equation}

\begin{equation}
\tag{2.b} A^{\ast} p+\mu^{\ast} = y^{\ast}-z\text{\quad in~}
H^{-1}\left(\Omega\right),
\end{equation}

\begin{equation}
\tag{3.b} \left\langle -\nu \Delta \varphi^{\ast} +\nu
\varphi^{\ast} + \mu^{\ast} ,\chi -\varphi^{\ast} \right\rangle \geq
0\text{\quad for all~}\chi \in \mathcal{U}_{ad},
\end{equation}

\begin{equation}
\tag{4.b} \mu^{\ast} \left( y^{\ast}-\varphi^{\ast} \right)
=0\text{\quad a.e. in~}\Omega,
\end{equation}

\begin{equation}
\tag{5.b} p^{\ast}\xi^{\ast} =0\text{\quad a.e. in~}\Omega,
\end{equation}

\begin{equation}
a^{*}\left(p^{\ast},p^{\ast}\right) -\left( z-y^{\ast },p^{\ast
}\right) \leq 0, \tag{6.a}
\end{equation}

\begin{equation}
\tag{7.b}
\left\langle \mu^{\ast} ,p^{\ast}\phi \right\rangle\geq 0\text{\quad for all~}%
\phi \in W^{1,\overline{q}}\left( \Omega \right), \text{\quad
with~}\phi \geq 0,\text{~and~}\overline{q}>n.
\end{equation}

\end{subequations}

Indeed, they studied the following optimal control problem
$(\widetilde{\mathcal{P}})$

\begin{equation*}
\left\langle Ay-f,\phi -y\right\rangle \geq 0\text{\quad for all
~}\phi \in \mathcal{K}\left(\varphi\right),
\end{equation*}%

with

\begin{equation*}
f\in \Ld,\quad \varphi \in \mathcal{U}_{ad}\text{\quad with~}
\varphi \leq 0 \text{\quad in ~}\partial \Omega,
\end{equation*}%

such that the cost functional is given by

\begin{equation*}
\widetilde{J}(\varphi
)\mathrel{\mathop:}=\tfrac{1}{2}\int\nolimits_{\Omega }\left(
\mathcal{T}\left( \varphi \right) -z\right) ^{2}d{
x}+\tfrac{\nu}{2}\left( \int\nolimits_{\Omega }\left(
|\varphi|^{2}+ |\nabla \varphi|^{2} \right)d{x}\right),
\end{equation*}

We notices a likeness between the two systems $(\mathcal{S})$ and
$(\widetilde{\mathcal{S}})$ , exepted for the equations (3.a) and
(3.b) (are respectively the differential of the objective function
$J$ and $\widetilde{J}$); i.e. in \cite{Kunish}, the authors
treated the optimal control problem $(\widetilde{\mathcal{P}})$,
such that $\varphi$ belong to $\mathcal{U}_{ad}
\mathrel{\mathop:}= \{ \varphi \in X: \varphi\left(x\right) \geq
0, \mbox{ on }
\partial\Omega \mbox{ and } -a\left(\varphi,v\right)+\left(f,v\right)
\leq \left(\overline{\lambda},v\right) \mbox{ for all } v \in V
\mbox{ with } v \geq 0\}$ with $H^{1}$-obstacle (where
$\overline{\lambda} \in \Ld \mbox{ satisfying } \overline{\lambda}
\geq 0 $  a.e. on $\Omega$), and in our work, we had study the
optimal control problem $(\mathcal{P})$, where $\varphi$ is in
$\mathcal{U}_{ad} \mathrel{\mathop:}= H^{2}\left(\Omega\right)\cap
H_{0}^{1}\left(\Omega\right)$, with $H^{2}$-obstacle.
\end{rem}

%---------------------------------------------------------------------------
\section*{Conclusion}
%---------------------------------------------------------------------------

In this work, we treated the theoretical aspect of the problem
$(\mathcal{P})$, we proved the existence of optimal solutions and
constructed a necessary optimality conditions system. Additional
optimal obstacle regularity has been also provided. Currently we
study the numerical aspect of the problem $(\mathcal{P})$, via a
numerical strategy based on the direct resolution of the
optimality system and using a fixed point algorithm.

\thanks{The author is grateful to Prof. M. Bergounioux for their
instructive suggestions.}

\end{document}